\documentclass{article}

\newtheorem{fed}{\textbf{Definition}}[section]
\newtheorem{thm}[fed]{\textbf{Theorem}}
\newtheorem{lemma}[fed]{\textbf{Lemma}}

\newtheorem{rem}[fed]{\textbf{Remark}}
\newtheorem{prop}[fed]{\textbf{Proposition}}
\newtheorem{cor}[fed]{\textbf{Corollary}}
\usepackage{amssymb,bbm,graphicx,epsfig,psfrag,epic,eepic,latexsym}
\usepackage{amsmath}
\usepackage{mathrsfs}
\usepackage
{color}

\newcommand{\R}{\mathbb{R}}

\begin{document}
\title{Existence of either a periodic collisional orbit or infinitely many consecutive collision orbits in the planar circular restricted three-body problem}
\author{Urs Frauenfelder, Lei Zhao}

\maketitle

\begin{abstract} In the restricted three-body problem, consecutive collision orbits are those orbits which start and end at collisions with one of the primaries. Interests for such orbits arise not only from mathematics but also from various engineering problems. In this article, using Floer homology, we show that there are either a periodic collisional orbit, or infinitely many consecutive collision orbits in the planar circular restricited three-body problem on each bounded component of the energy hypersurface for Jacobi energy below the first critical value.
\end{abstract}


\section{Introduction}

In this note we explain that below the first critical value for any mass ratio in each bounded component of the circular restricted three body problem there are infinitely many consecutive collision orbits. 

We first give some motivation for this question. 
In his lectures ``Hamiltonian dynamics and symplectic rigidity'' held at the workshop ``$J$-holomorphic Curves
in Symplectic Geometry, Topology and Dynamics'' at the CRM in Montr\'eal 2013, Helmut Hofer strongly emphazised the importance of finding symplectic relevant sets in order to solve the old dynamical problem of travelling from $A$ to $B$. In \cite{albers-frauenfelder-koert-paternain}, it was {shown} that below the first critical value, the {Moser-regularized} bounded components of the planar circular restricted three body problem can be interpreted as fiberwise starshaped hypersurfaces in the cotangent bundle over the two dimensional sphere. Moser regularization \cite{moser} first interchanges the role of fiber and base, i.e., momentum is interpreted as position and position as momentum, and then adds a fiber at infinity which corresponds to collisions where the momentum explodes. {In symplectic geometry, Lagrangian submanifolds are important submanifolds of symplectic manifolds in this sense.} A natural class
of Lagrangians in a cotangent bundle are the fibers. Since Moser regularization interchanges the role of position and momentum, fibers have a {slightly} strange meaning in the original coordinates, namely they correspond to a fixed momentum at an arbitrary location. Chords from one fibre to another one correspond to the problem of starting with a given momentum at an arbitrary place and ending up with another given momentum at an arbitrary place.  
We do not know if this kind of issue {arises anywhere} in engineering problems.

However, if we take our two base points both as the point at infinity, {such a chord then} corresponds to a consecutive collision orbit, which starts and ends in collisions. {At a first glance}, one might think that such a problem is only of interest to people tired of life. Nevertheless this problem is of use in engineering problems, see for example \cite{santos-prado}. {Indeed, with a slight perturbation, one obtains from such a consecutive collision orbit an orbit which slightly avoids collisions, which} are of broad use in engineering problems. A first application is the Oberth maneuver \cite{oberth} which is basically a very clever application of the binomial formula. Indeed, comparing the kinetic energy before and after a burn $\Delta v$ one obtains
$$\tfrac{1}{2}(v+\Delta v)^2-\tfrac{1}{2} v^2=\tfrac{1}{2} \Delta v^2+ v \cdot \Delta v.$$
The term of interest is the mixed term $v \cdot \Delta v$. It shows that if the velocity is high with a small
$\Delta v$-burn one can gain an enormous amount of energy. Different applications are gravitational slingshots which were for example used in the Voyager missions, see \cite[Chapter 8.9]{curtis}.

As has been explained by M. H\'enon \cite{henon}, such orbits are useful as well for space mission designs with a space explorer: with such an orbit, it is direct to reclaim the space explorer, or, with a slight modification of this orbit so that it does not collide but stay very close to the Earth, the recorded informations of the observation maybe transmitted back to the Earth with radio signals in the best condition. The interest of these orbits is more evident for spaceships, where as one can always command a spaceship on any orbit to turn back, there is always a risk of failure for tele-commanding, and the consummation of fuel is inevitable for such a turning-back, which could also be a problem some times. These orbits are thus optimal in this perspective also for these spaceships.

M. H\'enon made a detailed analysis on the existence of such orbits in the limiting case of a zero-mass Earth and a unit-mass Sun in a fixed reference frame. Finally, H\'enon has remarked that such orbits can be interpreted as the limiting case of the periodic orbits of second species of Poincar\'e in the restricted three-body problem.

\section{Consecutive collision orbits}
We consider the planar restricted circular three-body problem, in a proper rotating frame with a further shift of one of the primaries to locate at the origin $O$. We denote by $E:=(1,0) \in \R^{2}$ the position of the other primary. The masses of the primary are $1-\mu, \, \mu$ respectively with $0<\mu<1$. The Hamiltonian of the system is given by
$$H(p, q)=\dfrac{\|p\|^{2}}{2} + (p_{1} q_{2}- p_{2} (q_{1}-\mu)) + \dfrac{\mu}{\|q-E\|}+\dfrac{1-\mu}{\|q\|},$$
where $(p, q) \in \R^{2} \times (\R^{2} \setminus \{O, E\})$.

\begin{fed}  An orbit $x(t): (0, \tau) \subset \R \to \R^{2} \setminus \{O, E\}$  is called \emph{consecutive collision orbit} if it satisfies
$$\lim_{t \to 0} x(t)=\lim_{t \to \tau} x(t)=O.$$ 
\end{fed}
This simply means that the orbit goes from a collision with $O$ to another collision with $O$, whence the name. {Note that since our system is autonomous, we have put in our definition starting time instant $0$ to avoid multiple counting of time shifts of one orbit.}
\\ \\
{We know that in our system, two body collisions can be regularized. Therefore an orbit ending in a collision can be continued over the collision as a collision-ejection orbit. It might even happen that after a collision it collides again and gets continued further. This even does not rule out the possibility that one gets in this way a periodic orbit.
We refer to a periodic orbit containing at least one collision with $O$ as a \emph{periodic collision orbit}.}

In this paper, we shall show that

{\begin{thm}\label{thm:1.2} For any $0<\mu <1$, there exist a periodic collision orbit or infinitely many consecutive collision orbits in the planar circular restricted three-body problem for any energy hypersurface, for Jacobi energy below the first critical value.
 \end{thm}}
{\begin{rem}
By interchanging the roles of $O$ and $E$ it follows that as well in the bounded component containing $E$
there is a periodic collision orbit or infinitely  many consecutive collision orbits.
\end{rem}}

{To prove the Theorem we shall realize such orbits as chords in phase space in the sense of contact geometry via Moser regularization. A crucial ingredient is the fact that one can interpret the Hamiltonian flow in phase space as a Reeb flow \cite{albers-frauenfelder-koert-paternain}, {and the collision set is transformed into a Legendrian submanifold}. {Our result thus follows from} the existence of infinitely many Reeb chords in this context{, which} is a well-known fact from contact geometry and will be recalled in the sequel to complete the proof. }

{In our setting, a Reeb chord connection {the Legendrian of collisions} corresponds to an orbit which
starts and ends in a collision, {and} might go in between to additional collisions. {Should a periodic collision orbit exist, then} by iterating this orbit orbit one {also} obtains infinitely many Reeb chords,{ and in this case, with our normalization, we get only a periodic collision orbit. However this is a degenerate case, since} the collision set is a one dimensional Legendrian submanifold in the three dimensional energy hypersurface in four dimensional phase space. {Since} periodic
orbits are one dimensional objects{,} they generically avoid {this} Legendrian. This means that after {a proper} arbitrarily small perturbation of the contact form defining the Reeb flow of the restricted three body problem there are no
periodic collision orbits left and one gets {in this case} infinitely many consecutive collision orbits. } {In our concrete systems, however, it seems to be a delicate issue to exclude the existence of periodic collision orbits. }

Our study is limited to the case where the Jacobi energy lies below the first critical value, so that the component of the (regularized) energy hypersurface containing the singularity $O$ does not contain any other singularities, thus is a priori irrelevant to the periodic orbits of the second species which wind around the other primary. On the other hand, we do not impose any conditions on the masses of the primaries, and this result can be also obtained for more general Stark-Zeeman systems with the same approach.

We note that the existence of a {continuum} of such orbits in the planar circular restricted three-body problem may also be obtained by first going back to the fixed reference frame in which the two primaries undergo circular motion. Such orbits may be obtained by minimizing the action functional with initial and end configurations corresponding to consecutive collisions, and with arbitrary time interval. By Marchal's theorem for the massless body as stated in \cite{marchal}, the minimizers will be collision-free and thus give rise to classical solutions of our system outside the initial and end configurations. Note that with this approach we do not get any information on the Jacobi energy of these orbits in the rotating frame.

\section{The Moser regularization}

The singularity of $H$ at $q=O$ can be regularized with the following method of Moser \cite{moser}. 

We first fix an energy hypersurface $$\Sigma_{f}:=\{H+f=0\}$$ and change time on this hypersurface by multiplying $H+f$ with $\|q\|$, thus obtain
\small$$K(p,q):=\|q\| (H +f)=\dfrac{\|q\|(\|p\|^{2}+1)}{2} + (p_{1} q_{2}- p_{2} (q_{1}-\mu)+f-1/2)\|q\| + \dfrac{\mu \|q\|}{\|q-1\|}+(1-\mu).$$\normalsize

The flow of $K$ in $\{K=0\}$ agrees with that of 
\small$$\check{K}:=|K-(1-\mu)|^{2}/2=\Bigl(\dfrac{\|q\|(\|p\|^{2}+1)}{2} + (p_{1} q_{2}- p_{2} (q_{1}-\mu)+f-1/2)\|q\| + \dfrac{\mu \|q\|}{\|q-1\|}\Bigr)^{2}/2$$\normalsize
in the energy hypersurface 
$$\{\check{K}=(1-\mu)^{2}/2\}$$ up to {constant reparametrization of time}.
We may now pull-back the function $\check{K}$ by the stereographic projection from the north pole $N$ of the sphere $\mathbb{S}^{2} \in \R^{3}$ to $\R^{2} \times \{0\} \in \R^{3}$ while regarding $p$ as the positions and $-q$ the conjugate momenta. It is seen that with the energy restricted to $(1-\mu)^{2}/2$, the pull-backed flow extends smoothly through the fibre {$T^*_{N} \mathbb{S}^{2}$} over $N$ in $T^{*}\mathbb{S}^{2}$. The completion of the flow thus gives a regularization of the flow of $H$ in $\Sigma_{f}$.

Let us denote the relevant connected component of this completed energy hypersurface in $T^{*}\mathbb{S}^{d}$ by $\check{\Sigma}$. The intersection {$\check{\Sigma} \cap T^*_{N} \mathbb{S}^{d}$} corresponds to all physical collisions $q=O$ within $\check{\Sigma}$. 

In \cite{albers-frauenfelder-koert-paternain}, it is shown that up to the first critical value of $H$, the hypersurface $\check{\Sigma}$ is fibrewise star-shaped {with respect to the zero-section} and thus admit a contact 1-form which is simply the restriction of the Liouville 1-form of $T^{*} \mathbb{S}^{d}$ to this hypersurface. It follows that {$\check{\Sigma} \cap T^*_{N} \mathbb{S}^{d}$ } is a Legendrian submanifold inside $\check{\Sigma}$. Up to time parametrization, a consecutive collision orbit is thus realized in $\check{\Sigma}$ as a Reeb orbit connecting this Legendrian submanifold with itself. The existence of infinitely many consecutive collision orbits thus follows from the existence of infinitely many Reeb chords. In the next section, we shall explain this result on Reeb chords.

\section{Existence of Reeb chords}
\subsection{Reeb chords}

Assume that $N$ is a closed connected manifold and 
$$\Sigma \subset T^* N$$
is a closed fiberwise starshaped hypersurface. Fiberwise starshaped means that the Liouville
vector field $p \partial_p$ is transverse to $\Sigma$ so that the restriction of the Liouville
one-form $\lambda=p dq$ to $\Sigma$ is a contact form on $\Sigma$. We further pick a base point
$q_0 \in N$. Then the cotangent fibre $T^*_{q_0} N$ is a Lagrangian submanifold of $T^* N$ with respect
to the standard symplectic form $\omega=d \lambda=dp \wedge dq$. Moreover, the restriction of
the Liouville one-form $\lambda$ to $T^*_{q_0} N$ vanishes. Further note that
$$\mathscr{L}:=\Sigma \cap T_{q_0}^* N \subset \Sigma$$
is a Legendrian submanifold. Abbreviate by
$$\lambda_\Sigma:=\lambda|_{\Sigma} \in \Omega^1(\Sigma)$$
the contact form on $\Sigma$ obtained by restriction of the Liouville one-form to $\Sigma$. Define
the Reeb vector field $R \in \Gamma(T\Sigma)$ by
$$\lambda_\Sigma(R)=1, \quad \iota_R d\lambda_\Sigma=0.$$
\begin{fed}
A \emph{Reeb chord} $(x,\tau) \in C^\infty([0,1],\Sigma) \times(0,\infty)$ from $\mathscr{L}$ to $\mathscr{L}$ is a solution of the problem
$$\left\{\begin{array}{cc}
\partial_t x(t)=\tau R(x(t)) & t \in [0,1],\\
x(0),\,\,x(1) \in \mathscr{L}. & 
\end{array}\right.$$
\end{fed}
If one reparametrizes a Reeb chord $(x,\tau)$ and defines
$$x_\tau \colon C^\infty([0,\tau],\Sigma), \quad t \mapsto x\big(\tfrac{t}{\tau}\big)$$
one obtains a solution of the problem
$$\left\{\begin{array}{cc}
\partial_t x_\tau(t)=R(x_\tau(t)) & t \in [0,\tau],\\
x_\tau(0),\,\,x_\tau(\tau) \in \mathscr{L}. & 
\end{array}\right.$$
In view of this fact we refer to $\tau$ as the \emph{period} of the Reeb chord. 
\\ \\
The following theorem is known to experts in Floer homology.

\begin{thm}\label{infchord}
There exist infinitely many Reeb chords from $\mathscr{L}$ to $\mathscr{L}$.
\end{thm}

The existence of infinitely many consecutive collision solutions follows directly from this theorem.

For the readers' convenience we explain a proof of this result in terms of Rabinowitz Floer homology.
Alternative proofs could be provided by the use of wrapped Floer homology. 
\subsection{Rabinowitz Floer homology}

We first explain how Reeb chords from $\mathscr{L}$ to $\mathscr{L}$ can be interpreted variationally
as critical points of Rabinowitz action functional. For a smooth function $H \in C^\infty(T^*N, \mathbb{R})$
the Hamiltonian vector field $X_H \in \Gamma(T T^*N)$ is defined implicitly by the condition
$$dH= \omega( \cdot, X_H).$$
We pick a smooth function $H$ on the cotangent bundle of $N$ meeting the following conditions
\begin{description}
 \item[(i)] The differential $dH$ has compact support, i.e., $H$ is locally constant outside a compact subset of
$T^*N$.
 \item[(ii)] $\Sigma=H^{-1}(0)$, i.e., $\Sigma$ is the level set of $H$ to the value $0$.
 \item[(iii)] $X_H|_\Sigma=R$, i.e., the restriction of the Hamiltonian vector field of $H$ coincides
with the Reeb vector field on $\Sigma$. In particular, in combination with (ii) it follows that $0$ is a regular
value of $H$. 
\end{description}
We mention, that although the chain complex for Rabinowitz Floer homology depends on the choice of
the Hamiltonian $H$, the homology is independent of this choice and just depends on $\Sigma$ and
$\mathscr{L}$. Abbreviate
$$\mathcal{P}:=\mathcal{P}_{q_0}:=\big\{x \in C^\infty([0,1],T^*N): \,\,x(0),\,x(1) 
\in T_{q_0}^* N\big\}$$
the space of paths in $T^*N$ which start and end in the cotangent fibre over $q_0 \in N$. Rabinowitz action functional is defined as
$$\mathcal{A}^H \colon \mathcal{P} \times (0,\infty) \to \mathbb{R}, \quad
(x,\tau) \mapsto \int_0^1 x^* \lambda-\tau \int_0^1 H(x(t)) dt.$$
The first term is just the area functional and $\tau$ can be interpreted as a Lagrange multiplier, so that 
Rabinowitz action functional is the Lagrange multiplier functional of the area functional to the constraint
given by the vanishing of the mean value of $H$. Critical points of $\mathcal{A}^H$ are solutions
$(x,\tau) \in \mathcal{P} \times (0,\infty)$ of the problem
\begin{equation}\label{rab1}
\left\{\begin{array}{cc}
\partial_t x(t)=\tau X_H(x(t)) & t \in [0,1],\\
\int_0^1 H(x(t))dt=0. & 
\end{array}\right.
\end{equation}
The first equation is obtained by differentiating $\mathcal{A}^H$ with respect to the first variable, namely the path, while the second equation one gets by differentiating Rabinowitz action functional with respect to the second variable, namely the Lagrange multiplier. However note, that in view of the first equation by preservation of energy the Hamiltonian $H$ is constant along $x$. Combining this fact with the second equation we see that
solutions of problem (\ref{rab1}) are in one-to-one correspondence with solutions of the following problem
\begin{equation}\label{rab2}
\left\{\begin{array}{cc}
\partial_t x(t)=\tau X_H(x(t)) & t \in [0,1],\\
H(x(t))=0, &  t \in [0,1],
\end{array}\right.
\end{equation}
i.e., the mean value constraint is equivalent to a pointwise constraint. However, in view of properties (i) and
(ii) of $H$ we observe, that solutions of (\ref{rab2}) are precisely Reeb chords from $\mathscr{L}$ to
$\mathscr{L}$. In particular, we have proved the following proposition
\begin{prop}\label{crit}
Critical points of $\mathcal{A}^H$ are in one-to-one correspondence with Reeb chords from 
$\mathscr{L}$ to $\mathscr{L}$.
\end{prop}
Given a critical point $(x,\tau)$ of $\mathcal{A}^H$, i.e., a Reeb chord from $\mathscr{L}$ to
$\mathscr{L}$ we compute its action value as
$$\mathcal{A}^H(x,\tau)=\int_0^1 x^* \lambda=\int_0^1 \tau \lambda_\Sigma(R(x(t)))dt=\tau.$$
Hence we obtain the following proposition
\begin{prop}\label{act}
The action value of Rabinowitz action functional at a critical point $(x,\tau)$ is $\tau$, namely the period
of the Reeb chord. 
\end{prop}
To construct the chain complex of Rabinowitz Floer homology we need an additional assumption on 
$\Sigma$. We denote for $t \in \mathbb{R}$ by $\phi^t_R \colon \Sigma \to \Sigma$ the flow of
the Reeb vector field on $\Sigma$.
\begin{fed}
$\Sigma$ is called \emph{regular} if for every Reeb chord $(x,\tau)$ from $\mathscr{L}$ to $\mathscr{L}$
it holds that
$$d \phi^\tau_R(x(0)) T_{x(0)} \mathscr{L} \cap T_{x(1)}\mathscr{L}=\{0\}.$$
\end{fed}
If $\Sigma$ is regular, the kernel of the Hessian of Rabinowitz action functional is trivial at every critical point, and hence Rabinowitz action functional is Morse. In particular, its critical points are isolated. We define a graded 
vector space
$$CM_*(\mathcal{A}^H):=\mathrm{crit}(\mathcal{A}^H) \otimes \mathbb{Z}_2,$$
i.e., the $\mathbb{Z}_2$-vector space generated by critical points of Rabinowitz action functional. The grading
is obtained by the transversal Maslov index at a chord. One defines a boundary operator
$$\partial \colon CM_*(\mathcal{A}^H) \to CM_{*-1}(\mathcal{A}^H)$$
by counting the number of gradient flow lines modulo two between critical points of index difference one as in
finite dimensional Morse homology. Here the gradient of Rabinowitz action functional is defined as in classical Floer theory via a metric obtained from a family of $\omega$-compatible almost complex structures on
$T^*N$. We define Rabinowitz Floer homology as
$$RFH_*^+(\Sigma,\mathscr{L}):=HM_*(\mathcal{A}^H):=\frac{\mathrm{ker} \partial}{\mathrm{im} \partial}.$$
As the notation suggests the resulting homology is independent of the choice of the Hamiltonian $H$ as well
as on the choice of the family of $\omega$-compatible almost complex structures needed to define the gradient \cite{cieliebak-frauenfelder, merry}. The superscript $+$ is added, because we suppose that $\tau$ only takes positive values. In the full Rabinowitz Floer homology $RFH_*(\Sigma,\mathscr{L})$ the Lagrange multiplier $\tau$ is allowed to assume
any real value. In the full case critical points with $\tau=0$ correspond to constant chords, i.e., points on $\mathscr{L}$,
and critical points with $\tau<0$ correspond to chords traversed backwards. 
\\ \\
Abbreviate 
$$\Omega_{q_0}:=\big\{q \in C^\infty([0,1],N): q(0)=q(1)=q_0\big\}$$
the based loop space of $N$. Note that $q_0$ interpreted as a constant loop becomes itself an element
of $\Omega_{q_0}$. We abbreviate by $H_*(\Omega_{q_0},q_0)$ the homology of the based loop space relative to the constant loop $q_0$ with $\mathbb{Z}_2$-coefficients. The following result is due to Merry \cite{merry}.
\begin{thm}[Merry]\label{merthe}
Assume that $\Sigma$ is regular, then
$$RFH_*^+(\Sigma,\mathscr{L})=H_*(\Omega_{q_0},q_0),$$
i.e., there is a canonical isomorphism
$$\zeta_\Sigma \colon RFH_*^+(\Sigma,\mathscr{L}) \to H_*(\Omega_{q_0},q_0).$$
\end{thm}
A Corollary of Merry's theorem is the following.
\begin{cor}\label{mercor}
Assume that $\Sigma$ is regular, then
$$\mathrm{dim} (RFH_*^+(\Sigma,\mathscr{L}))=\infty.$$
\end{cor}
\textbf{Proof: }In the case that the fundamental group $\pi_1(N)$ is infinite, the based loop space has
in view of the equality $\pi_1(N)=\pi_0(\Omega_{q_0})$
infinitely many connected components and therefore already the zeroth Betti number
satisfies $b_0(\Omega_{q_0},q_0)=\infty$.

If $N$ is simply connected, then $\dim(H_*(\Omega_{q_0},q_0))=\infty$ follows from Serre's spectral sequence. Finally if $N$ has finite fundamental group, its universal cover $\widetilde{N}$ is still closed. Moreover,
the based loop space of $N$ decomposes into connected components indexed by the fundamental group each one homotopic to the based loop space of the universal cover $\widetilde{N}$, and the result follows again from Serre's spectral sequence. \hfill $\square$

To obtain infinitely many consecutive collision orbits for the planar case one just needs the case $N=S^2$, where the $\mathbb{Z}_2$-Betti numbers satisfy
$$b_k(\Omega_{q_0},q_0)=1, \quad k \in \mathbb{N}.$$

In view of Proposition~\ref{crit} and Corollary~\ref{mercor} in the case that $\Sigma$ is regular Theorem~\ref{infchord} follows now immediately from the Morse inequalities for Rabinowitz Floer homology. 
To prove the theorem as well in the case where $\Sigma$ is not necessarily regular, so that Rabinowitz action functional does not need to be Morse and Rabinowitz Floer homology cannot be defined directly, we introduce spectral invariants. 
{\subsection{Spectral invariants}}
We first define spectral invariants under the hypothesis that $\Sigma \subset T^*N$ is a regular closed fiberwise starshaped hypersurface and then get rid in a second step of the regularity assumption. If $\Sigma$
is regular, Rabinowitz action functional $\mathcal{A}^H$ is Morse and for $k \in \mathbb{Z}$ a vector
$\xi \in CM_k(\mathcal{A}^H)$ is a sum
$$\xi=\sum_{\substack{(x,\tau) \in \mathrm{crit}(\mathcal{A}^H)\\
\mu(x,\tau)=k}}\xi_{(x,\tau)}(x,\tau)$$
where $\mu(x,\tau)$ is the Maslov index of the Reeb chord $(x,\tau)$. The coefficients $\xi_{(x,\tau)}$ belong to the field $\mathbb{Z}_2$ and only finitely many coefficients are $1$. We set
$$\sigma(\xi):=\max\{\tau: \xi_{(x,\tau)}=1\}.$$
By Proposition~\ref{act} this corresponds to the maximal action value of $\mathcal{A}^H$ on the formal
sum of its critical points contributing to $\xi$. If $\alpha \in RFH^+(\Sigma,\mathscr{L})$ we set
$$\sigma(\alpha):=\min\big\{\sigma(\xi): \alpha=[\xi]\big\}.$$
Abbreviate
$$\mathscr{S}:=\big\{\Sigma \subset T^*N: \Sigma\,\,\textrm{closed and fiberwise starshaped}
\big\}$$
and 
$$\mathscr{S}_{\mathrm{reg}}:=\big\{\Sigma \in \mathscr{S}: \Sigma\,\,\textrm{regular}\big\}.$$
We endow $\mathscr{S}$ with the $C^0$-topology. By Sard's theorem it holds that
$$\mathscr{S}_{\mathrm{reg}} \subset \mathscr{S}$$
is dense. {We fix an element} 
$$\beta \in H_*(\Omega_{q_0}, q_0).$$
{For} $\Sigma \in \mathscr{S}_{\mathrm{reg}}${,} let
$$\zeta_\Sigma \colon RFH^+(\Sigma,\mathscr{L}) \to H_*(\Omega_{q_0},q_0)$$
be the isomorphism from Theorem~\ref{merthe}. 
We define
$$\rho_\beta \colon \mathscr{S}_{\mathrm{reg}} \to (0,\infty)$$
for $\Sigma \in \mathscr{S}_{\mathrm{reg}}$ by
$$\rho_\beta(\Sigma)=\sigma(\zeta_\Sigma^{-1}(\beta)).$$
One can show \cite{albers-frauenfelder} that the function $\rho_\beta$ is locally Lipschitz continuous.
Because $\mathscr{S}_{\mathrm{reg}} \subset \mathscr{S}$ is dense there exists a unique continuous
extension 
$$\overline{\rho}_\beta \colon \mathscr{S} \to (0, \infty),$$
i.e., $\overline{\rho}_\beta$ is characterized by the properties
\begin{description}
 \item[(i)] $\overline{\rho}_\beta|_{\mathscr{S}_{\mathrm{reg}}}=\rho_\beta,$
 \item[(ii)] $\overline{\rho}_\beta$ is continuous. 
\end{description}
For $\Sigma \in \mathscr{S}$ define the \emph{spectrum} of $\Sigma$ by
$$\mathfrak{S}(\Sigma)=\big\{\tau: (x,\tau)\,\,\textrm{Reeb chord from}\,\mathscr{L}\,\textrm{to}\,
\mathscr{L}\big\}.$$
Note that by Proposition~\ref{act} the spectrum of $\Sigma$ coincides with the action spectrum of {the}
Rabinowitz action functional. {We refer} 
$$\mathfrak{S}:=\{(\Sigma,\tau): \Sigma \in \mathscr{S},\,\,\tau \in \mathfrak{S}(\Sigma)\big\} \subset \mathscr{S} \times (0,\infty)$$
as the \emph{spectral bundle} over $\mathscr{S}${, with the} canonical projection
$$\pi \colon \mathfrak{S} \to \mathscr{S}, \quad (\Sigma,\tau) \mapsto \Sigma.$$
The following lemma justifies {our terminology \emph{spectral invariant}} for $\overline{\rho}_\beta$. 
\begin{lemma} \label{speclem}
The map $\overline{\rho}_\beta$ is a section from $\mathscr{S}$ to $\Sigma$, i.e., 
\begin{equation}\label{spec}
\overline{\rho}_\beta(\Sigma) \in \mathfrak{S}(\Sigma), \quad \forall\,\,\Sigma \in
\mathscr{S}.
\end{equation}
\end{lemma}
\textbf{Proof: } If $\Sigma \in \mathscr{S}_{\mathrm{reg}}$, then $\overline{\rho}_\beta(\Sigma)=\rho_\beta(\Sigma)$ and (\ref{spec}) follows from the definition of $\rho_\beta$. 
For the general case, there exists a sequence $\Sigma_\nu \in \mathscr{S}_{\mathrm{reg}}$,
for $\nu \in \mathbb{N}$, which converges to $\Sigma$. Because $\Sigma_\nu \in \mathscr{S}_{\mathrm{reg}}$ we have
$$\overline{\rho}_\beta(\Sigma_\nu)=\rho_\beta(\Sigma_\nu) \in
\mathfrak{S}(\Sigma_\nu), \quad \nu \in \mathbb{N}.$$
In particular, there exists a sequence of Reeb chords $(x_\nu,\tau_\nu)$ for $\Sigma_\nu$ such that
$$\overline{\rho}_\beta(\Sigma_\nu)=\tau_\nu.$$
Because $\overline{\rho}_\beta$ is continuous and $\Sigma_\nu$ converges, we conclude that
the sequence $\tau_\nu$ is bounded. Therefore by the Theorem of Arzela-Ascoli there exists a subsequence
$\nu_j$ of $\nu$ and a Reeb chord $(x,\tau)$ of $\Sigma$ such that
$$\lim_{j \to \infty}(x_{\nu_j},\tau_{\nu_j})=(x,\tau).$$
Again by continuity of $\overline{\rho}_\beta$ we obtain
$$\tau=\lim_{j \to \infty}\tau_{\nu_j}=\lim_{j \to \infty} \overline{\rho}_\beta(\Sigma_{\nu_j})
=\overline{\rho}_\beta(\Sigma).$$
We have shown that  $\overline{\rho}_\beta(\Sigma) \in \mathfrak{S}(\Sigma)$. This finishes the proof of the Lemma. \hfill $\square$
\\ \\
If $\Sigma \in \mathscr{S}$ we define the \emph{homological spectrum} by
$$\mathfrak{H}(\Sigma)=\big\{\overline{\rho}_\beta(\Sigma): \beta \in H_*(\Omega_{q_0},q_0)\big\}.$$
By Lemma~\ref{speclem} we have
\begin{equation}\label{homspec}
\mathfrak{H}(\Sigma) \subset \mathfrak{S}(\Sigma).
\end{equation}
The following result was proved by Kang \cite{kang} using Corollary~\ref{mercor}.
\begin{thm}[Kang]\label{jungsoo}
For every $\Sigma \in \mathfrak{S}$ its homological spectrum $\mathfrak{H}(\Sigma)$ is unbounded.
\end{thm}
We are now in position to prove existence of infinitely many Reeb chords.
\\ \\
\textbf{Proof of Theorem~\ref{infchord}: }By Theorem~\ref{jungsoo} the homological spectrum of
$\Sigma$ is unbounded. By (\ref{homspec}) it follows that the spectrum of $\Sigma$ is unbounded.
In particular, there are infinitely many Reeb chords. \hfill $\square$
\\ \\
\textbf{Proof of Theorem~\ref{thm:1.2}:} In view of Theorem~\ref{infchord} it remains to show that the existence of infinitely many chords in the absence of a periodic collision orbit implies the existence of infinitely many consecutive collision orbits. Because the orbit of a first order ODE is uniquely determined by its initial condition in the absence of periodic
collision orbits the same consecutive collision orbit can occur at most one in a given orbit. Therefore there
have to exist infinitely many consecutive collision orbits in order to guarantee the existence of infinitely many chords. \hfill $\square$

\section{Some Remarks}
We {finish by} remarks.

\begin{rem} The existence of consecutive collision orbits is manifested in the rotating Kepler problem, as a limiting case of the restricted three-body problem, as all collisional orbits give rise to consecutive collision orbits. With our current approach, the result holds also true for Hill's problem, thanks to the result of J. Lee \cite{lee} on fiberwise convexity (and thus star-shapedness) {with respect to the zero section} of its regularization.
\end{rem}

\begin{rem} The planar restricted three-body problem is a special case in the family of Stark-Zeeman systems \cite{cieliebak-frauenfelder-koert}. Similar results may be obtained also for other Stark-Zeeman systems, as long as the fibrewise starshapedness {with respect to the zero section} of the corresponding component of the regularized energy hypersurface can be established.
\end{rem}

\begin{rem} A practical way to find such orbits is to use the shooting method of Birkhoff \cite{birkhoff}, by properly choosing initial positions on the axis of the primaries and initial velocity perpendicular to this axis. If such an orbit encounters a collision, it will have to encounter a collision in the past by the reflection symmetry with respect to the axis.
\end{rem}




\noindent
Urs Frauenfelder, University of Augsburg: {urs.frauenfelder@math.uni-augsburg.de} \\
Lei Zhao, University of Augsburg: l.zhao@math.uni-augsburg.de

\end{document}